\documentclass[oneside,english]{amsart}
\usepackage[latin1]{inputenc}
\usepackage{amssymb}
\usepackage{babel}

\makeatletter

\newtheorem*{thm*}{Theorem}

\title{On Multiple and Polynomial Recurrent extensions of infinite measure preserving transformations
}

\author{Tom Meyerovitch}
\makeatother

\newcommand{\ZD}{\mathbb{Z}^d}

\newtheorem{thm}{Theorem}[section]

\begin{document}
 \begin{abstract}We prove that multiple-recurrence  and polynomial-recurrence of invertible infinite measure preserving transformations
are both properties which pass to extensions.

%\noindent Resum\'{e}: Nous montrons que la r\'{e}currence  multiple
%et la r\'{e}currence  polyn\^{o}mial des transformations inversibles
%qui preservent une measure infinite sont deux propri\'{e}t\'{e}s qui
%passent aux extensions.

\end{abstract}

\maketitle

\section{Introduction and statement of result} A non-singular
transformation $T:X \to X$ of  a measure space $(X,\mathcal{B},\mu)$
is called \emph{$d$-recurrent} ($d \in \mathbb{N}$) if for any $A
\in \mathcal{B}$ of positive measure, there exists an integer $k \ge
1$ such that
\begin{equation}\label{eq:mult-rec} \mu(\bigcap_{i=0}^{d}T^{-ik}A)>0
\end{equation}
$T$ is \emph{multiply recurrent} if it is $d$-recurrent for all $d
\ge 1$. %We call a set $ A \in B$ for which equation
%\eqref{eq:mult-rec} holds  a $k$-recurrent set.

\noindent $T$ is called \emph{polynomially-recurrent} if for any $d
\ge 1$ and
 any polynomials $p_1,\ldots,p_d \in \mathbb{Z}[x]$ such that
$p_i(0)=0$ for all $1\le i \le d$, and any $A \in \mathcal{B}$ of
positive measure, there exists an integer $k  \ne 0$ such that
\begin{equation}\label{eq:poly-rec} \mu(\bigcap_{i=0}^{d}T^{-p_i(k)}A)>0
\end{equation}

We say that $T:X \to X$ is an extension of a measure preserving
transformation $S:Y \to Y$ (of the measure space
$(Y,\mathcal{C},\nu)$) if there is a measurable $\pi:X \to Y$ such
that $\pi\circ T = S \circ \pi$ and $\mu \circ \pi^{-1} = \nu$.

Furstenburg \cite{furst_mr_sz} gave an ergodic-theoretical proof
that any finite-measure preserving system is multiply recurrent,
giving an alternative proof of Szemer{\'e}di's theorem
\cite{szem_arithm}. Bergelson and Leibman proved that any
finite-measure preserving system is multiply recurrent
\cite{bergelson_leibman_96}, and deduced a theorem about existence
of polynomial configurations in subsets of positive density.

% previous results about multiple and polynomial recurrence of
%infinite measure preserving transformation transformations :

Combinatorial results about arithmetical progressions and polynomial
configurations in zero density sets (for example, the primes) are
generally more difficult to obtain using methods of classical
ergodic theory. There are old speculations about the relevance of
infinite-measure ergodic theory for such problems. In this
direction, Aaronson and Nakada \cite{aaro_naka_mr} formulated a
conjecture on infinite-measure preserving transformations which
holds assuming a positive solutions to a long standing conjuncture
of  Erd\"{o}s.

Several authors have studied multiple and polynomial properties of
infinite-measure preserving transformations: Eigen-Hajian-Halverson
\cite{eh98} constructed for each $d>0$ an ergodic infinite
measure-preserving transformation is that is $d$-recurrent but not
$(d+1)$-recurrent. Aaronson and Nakada \cite{aaro_naka_mr} give
necessary and sufficient conditions for $d$-recurrence of
Markov-Shifts.  Danilenko and Silva \cite{danilenko_silva_2004}
constructed measure preserving group actions with various multiple
and polynomial recurrence properties. In particular, there exist
infinite-measure preserving transformations which are polynomially
recurrent, and also measure preserving transformations which are
multiply-recurrent and not polynomially recurrent.

We prove the following results:
\begin{thm}
\label{thm:mult} If an invertible measure preserving transformation
$S$ is multiply recurrent, so is any measure preserving extension
$T$ of $S$.
\end{thm}

\begin{thm}\label{thm:poly}
If an invertible measure preserving transformation $S$ is
polynomially recurrent, so is any measure preserving extension $T$
of $S$.
\end{thm}
Theorem \ref{thm:mult} answers a question raised by Aaronson and
Nakada \cite{aaro_naka_mr}. A partial result was previously obtained
by Inoue \cite{inou_mr_ext} for isometric extensions. It is worth
noting that Inoue's result on multiple-recurrence for isometric
extensions does not require that transformations involved be
invertible.

It unknown if for some $k \ge 2$ there exist $n_k$ such that any
extension of any $n_k$-recurrent transformation is $k$-recurrent.

We remark that these results generalize essentially without
modification to measure preserving $\ZD$ actions, using the
multidimensional Szemer{\'e}di theorem and the Bergelson-Leibman
multidimensional polynomial-recurrence theorem.

Using an argument similar to the proof of our theorem
\ref{thm:mult},Furstenberg and Glasner obtained a
Szemered{'e}di-type theorem for `$\mathit{SL}(2,\mathbb{R})$
``$m$-stationary systems''. This result is described in
\cite{furst_glasner_msys}, and is based on a certain structure
theory of ``$m$-systems'', plus a specific ``multiple-recurrence
property'' related to $\mathit{SL}(2,\mathbb{R})$. The results of
this chapter were obtained independently of
\cite{furst_glasner_msys}.

\textbf{This work is a part of the Author's PhD.
 written under the supervision of Professor Jon Aaronson at Tel Aviv University.}
\section{Proof of theorem \ref{thm:mult}}
Recall the following formulation of Szemer{\'e}di's theorem:
\begin{thm*}\textbf{(Szemer{\'e}di's theorem - finitary version)}
\label{zm_finite_th} Let $l \in \mathbb{N}$ and $\delta>0$. For any
sufficiently large $L \in \mathbb{N}$, any $E \subset
\{1,\ldots,L\}$ with $|E| > \delta L$ contains a non-trivial
$l$-term arithmetic progression.
%There exists a function $L_{0}(\delta,l)$
%defined and finite for $\delta>0$ and a natural number $l$, so that
%if $L \ge L_0(\delta,l)$ and $S\subset \{1,2,\ldots,L\}$ is a subset
%with $|S| \ge \delta L$, then $S$ contains a non-trivial $l$-term
%arithmetic progression.
\end{thm*}
 Suppose $S:Y \to Y$ is an
invertible multiply recurrent measure preserving transformation, and
$T$ is an extension. We need to show that any set $A \in
\mathcal{B}$ with $0 < \mu(A) < \infty$ is multiply recurrent.

%The disintegration (``conditional expectation'') of $\mu$ over the
%$\sigma$-algebra $\pi^{-1}\C$ is a function $\rho:X \times \B \to
%\mathbb{R}_+$, such that for every $A \in \B$, the map $x \to
%\rho_x(A)$ is $C$-measurable, and $\rho$ satisfies the relation:
%$$\int_B \rho_x(A) d\mu(x)= \mu(A \cap B) ~~ \forall B\in
%\C , A \in \B$$
%Such a function is uniquely defined up to a
%$\mu$-null set in $X$.
%In case $\mu$ is a probability measure,
%$\rho_x(A) = E( 1_A | \C)(x)$ is the conditional expectation of $A$
%given the $\sigma$-algebra $\C$.

Denote by $\mu_y$ the conditional measure of $\mu$ given $y \in Y$.
This is $\nu$-almost everywhere defined by requiring that $y \to
\mu_y$ be $\mathcal{C}$-measurable and
$$\int_B \mu_y(A)d\nu(y)= \mu(A \cap \pi^{-1}B) ~~ \forall B\in \mathcal{C} , A \in \mathcal{B}$$

Since $\mu(T^{-1}A \cap T^{-1}B)=\mu(A \cap B)$ and
$T^{-1}\mathcal{C}=\mathcal{C}$, it follows that $\mu_{Sy}(A)=
\mu_y(T^{-1}A)$ for almost any $y \in Y$.

 Let $A \in \mathcal{B}$ with $0<\mu(A)$, and let
 $$B=B_\epsilon =\{y \in Y :\; \mu_y(A) >\epsilon \}$$
  We set some
$\epsilon>0$, so that $\mu(B_\epsilon)>0$. Note that $B \in
\mathcal{C}$, since $y \to \mu(A)$ is $\mathcal{C}$-measurable.

By multiple recurrence of $S$, for any $M\in \mathbb{N}$ there exist
$n\in \mathbb{N}$ such that $\nu(\bigcap_{j=0}^{M}S^{-jn}B)>0$.

 For
$y \in \bigcap_{j=0}^{M}S^{-jn}B$ and $0\le j \le M$, we have
$$\mu_y(T^{-jM}A)=\mu_{S^{jn}y}(A)>\epsilon.$$ Thus, for $y \in
\bigcap_{j=0}^{M}S^{-jn}B$, $$\int_X
\sum_{j=0}^{M}1_{T^{-jn}A}d\mu_y(x)> \epsilon M,$$ and so
$$ \int_{ \bigcap_{j=0}^{M}S^{-jn}B}\sum_{j=0}^{M}1_{T^{-jn}A}(x)d\nu(y) > \epsilon M \nu( \bigcap_{j=0}^{M}S^{-jn}B) $$
It follows   that there is a set $E \subset \{1,\ldots M\}$ with
$|E|
> \epsilon M$ with $$\mu(\bigcap_{j \in E}T^{-jn}A)>0.$$
 Choose
$M$ above large enough so that by  Szemer\'{e}di's theorem $E$
contains an arithmetic progression $(a+R,a+2R,\ldots,a+lR)$ of
length $l$. It follows that
$$\mu(\bigcap_{j=0}^{l-1}T^{-jR}A)>0,$$

and so  $T$ is multiply recurrent.

\section{Proof of theorem \ref{thm:poly}}

Our proof is based on the following theorem of V. Bergelson and A.
Leibman \cite{bergelson_leibman_96}, of which we state a finitary
version:
\begin{thm*}\textbf{(Bergelson-Leibman theorem, finitary version)}
Let
  $\{P_{i,j}(x)\}_{1\le i \le k, 1\le j\le d}$ be any polynomials with rational
  coefficients takeing on integer values in the integers and
  satisfying $P_{i,j}(0)=0$, and $\epsilon >0$.
  For any sufficiently large $N$, and any set
  $E \subset \{1,\ldots N\}^d$ with $|E| \ge \epsilon N^d$,
  there exists an integer $n$ and a vector $\overline{u} \in \ZD$ such that
  $\overline{u}
  + \sum_{j=1}^kP_{i,j}(x)e_d \in S$ for all $1\le i \le k$, where
  $\{e_1,\ldots,e_d\}$ are the standard basis of $\ZD$.
\end{thm*}

 Suppose
$S$ is an invertible polynomially-recurrent measure preserving
transformation, and $(X,\mathcal{B},\mu,T)$ is an extension. We need
to prove that for any $A \in \mathcal{B}$ with $\mu(A)>0$ and any
polynomials $p_1,\ldots,p_d \in \mathbb{Z}[x]$ with $p_i(0)=0$,
there exists $k \in \mathbb{Z}\setminus \{0\}$ such that equation
\eqref{eq:poly-rec} holds.

Write $p_i(x)= \sum_{j=1}^l a_{i,j}x^j$ For
$\overline{k}=(k_1,\ldots,k_l)$, let
$q_{\overline{k}}(x)=\sum_{j=1}^lk_jx^j$.

Find $B \in \mathcal{C}$ with $0< \nu(B) <+\infty$ and $\mu_y(A)
> \epsilon$ for all $y \in B$.

By polynomial-recurrence of $S$, there exists $n \in \mathbb{Z}$
such that

$$\nu(\bigcap_{k \in [N]^l}
S^{q_{\overline{k}}(n)}B)> 0,$$

with $N$ large enough so that the conclusion of the
Bergelson-Leibman theorem applies. Repeating the argument of the
previous proof, there exists a set $E \subset [N]^l$ of density at
least $\epsilon$ and a set $A' \subset A$ of positive measure such
that:
%$x\in T^{-p_{\overline{k}}(n)}A$ for all $\overline{k} \in E$. In
%other words,
$$\mu( \bigcap_{ \overline{k} \in E}T^{p_{\overline{k}}(n)}A') > 0$$

By the Bergelson-Leibman theorem, if $N$ is sufficiently large,
there exist $\overline{u} \in [N]^l$ and $r \in \mathbb{Z}$ such
that $\overline{u} + \sum_{j=1}^l a_{i,j}r^je_j \in E$ for every $1
\le i \le d$. It follows that

$$\mu(\bigcap_{i=1}^d T^{-p_i(rn)+U(n)}A)>0,$$
where $U(x)=\sum_{i=1}^lu_lx^l$, and so:
$$\mu(\bigcap_{i=1}^d T^{-p_i(rn)}A)>0$$

\bibliographystyle{abbrv}
\bibliography{ext_mult_rec}
\end{document}